\documentclass[12pt]{amsart}

\usepackage{amsmath,amsthm,amssymb}
\font\teneufm=eufm10
\font\seveneufm=eufm7
\font\fiveeufm=eufm5
\newfam\frakturfam

\textfont\frakturfam=\teneufm
\scriptfont\frakturfam=\seveneufm
\scriptscriptfont\frakturfam=\fiveeufm

\newtheorem{pr}{Proposition}
\newtheorem{ex}{Example}

\newtheorem{lemma}{Lemma}
\newtheorem{theorem}{Theorem}
\newtheorem{corol}{Corollary}

\def\bee{\begin{eqnarray}}
\def\bes{\begin{eqnarray*}}
\def\eee{\end{eqnarray}}
\def\ees{\end{eqnarray*}}

\def\a{\alpha}

\def\d{\partial}
\def\Proof{{\sl Proof.}\ }

\title{Automorphisms and derivations of a universal left-symmetric enveloping algebra}

\begin{document}
\date{}
\maketitle

\begin{center}
{\bf D. Zhangazinova, A. Naurazbekova, }\footnote{
L.N. Gumilyov Eurasian National University,
Astana, 010008, Kazakhstan,
e-mail: {\em naurazbekova82@gmail.com},}
{\bf U. Umirbaev}\footnote{
Department of Mathematics,
 Wayne State University,
Detroit, MI 48202, USA; Department of Mathematics, 
Institute of Mathematics and Mathematical Modeling, Almaty, 050010, Kazakhstan,
e-mail: {\em umirbaev@wayne.edu}}
\end{center}

\begin{abstract}
Let $A_n$ be an $n$-dimensional  algebra with zero multiplication over a field $K$ of characteristic $0$.  Then its universal (multiplicative) enveloping algebra $U_n$ in the variety of left-symmetric algebras is a homogeneous quadratic algebra generated by $2n$ elements $l_1,\ldots,l_n,r_1,\ldots,r_n$,  which contains both the polynomial algebra $L_n=K[l_1,\ldots,l_n]$ and the free associative algebra $R_n=K\langle r_1,\ldots,r_n\rangle$. We show that the automorphism groups of the polynomial algebra $L_n$ and the algebra $U_n$ are isomorphic for all $n\geq 2$, based on a detailed analysis of locally nilpotent derivations.  In contrast, we show that this isomorphism does not hold for $n=1$, and we provide a complete description of all automorphisms and locally nilpotent derivations of $U_1$.
\end{abstract}

\noindent


\tableofcontents

\section{Introduction}
\hspace*{\parindent}

Many papers in the theory of quantum groups are devoted to describing the automorphism groups of quadratic algebras \cite{AC92,AD96,CPWZ16,GL24,LL13,Yak13}. 
The Veronese subalgebras of free algebras \cite{AU23,AMU25} and the universal (multiplicative) enveloping algebras of finite-dimensional algebras \cite{Fleury,Joseph,NU23} provide a source of quadratic algebras with a rich group of automorphisms. For instance, the universal enveloping algebra of a finite-dimensional abelian Lie algebra is a polynomial algebra, and understanding its automorphisms is a central theme in affine algebraic geometry.

It is well known that automorphisms of the polynomial algebra $K[x,y]$ are tame \cite{Jung, Kulk}. Moreover, the automorphism group $\mathrm{Aut}(K[x,y])$
 of $K[x,y]$ admits an amalgamated free product structure \cite{Kulk, Shaf}, i.e.,
\bes
\mathrm{Aut}(K[x,y])=A\ast_{C}B,
\ees
where $A$ is the affine automorphism subgroup, $B$ is the triangular automorphism subgroup, and $C=A\cap B$. Similar results hold for free associative algebras \cite{Czer, M-L1} and free Poisson algebras in  characteristic zero \cite{M-L2}. Moreover, the automorphism groups of these algebras are isomorphic to the automorphism group of the polynomial algebra \cite{M-L2}. 

The automorphism groups of polynomial algebras \cite{Shestakov1, Shestakov2, Umirbaev1}, free associative algebras \cite{Umirbaev2, Umirbaev3}, and free Poisson algebras \cite{SZ24} in three variables over a field of characteristic zero cannot be generated by all elementary automorphisms, i.e., wild automorphisms exist.

In 1968, R. Rentschler \cite{Rentschler} proved that  locally nilpotent derivations of polynomial algebras in two variables over a field of characteristic zero are triangulable. An analogue of this result was proved for free associative algebras \cite{ANK,DM19,CS20} and for free Poisson algebras \cite{M-L2}. H. Bass \cite{Bass} provided an example of a non-triangulable locally nilpotent derivation of polynomial algebras in three variables over a field of characteristic zero.

An algebra $A$ over an arbitrary field $K$ with a bilinear product $x\cdot y$ is called a \textit{left-symmetric algebra} if the identity 
\begin{gather} \label{IN1}
(xy)z-x(yz)=(yx)z-y(xz)
\end{gather}
is satisfied for any $x,y,z\in A$. In other words, the associator $(x,y,z)=(xy)z-x(yz)$ is symmetric with respect to $x$ and $y$, i.e., $(x,y,z)=(y,x,z)$. The variety of left-symmetric algebras is \textit{Lie-admissible}, i.e., every left-symmetric algebra $A$ with respect to the operation $[x,y]=xy-yx$ is a Lie algebra. 

In 1994, D. Segal \cite{Segal} constructed a basis of  free left-symmetric algebras. In 2004, D. Kozybaev and U. Umirbaev \cite{Koz2} constructed a basis of the universal multiplicative enveloping algebra of a right-symmetric algebra. 

L. Makar-Limanov, D. Kozybaev, and U. Umirbaev \cite{Koz1} proved that all automorphisms of a free right-symmetric algebra of rank two are tame. Furthermore, A. Alimbaev, A. Naurazbekova, and D. Kozybaev \cite{ANK} demonstrated that the automorphism group of this algebra admits an amalgamated free product structure. They also proved that, in the case of a free right-symmetric algebra of rank two over a field of characteristic zero, every reductive automorphism group is linearizable, and every locally nilpotent derivation is triangulable.

In 1964, P. Cohn \cite{Cohn} proved that all automorphisms of finitely generated free Lie algebras over an arbitrary field are tame. In 1968, J. Lewin \cite{Lewin} extended this result to free algebras in Nielsen–Schreier varieties. Recall that a variety of algebras is called {\em Nielsen–Schreier} if every subalgebra of a free algebra in the variety is itself free. It is well known that the varieties of all nonassociative algebras \cite{Kurosh}, commutative and anticommutative algebras \cite{Shirshov1}, Lie algebras \cite{Shirshov2, Witt}, and Lie superalgebras \cite{Mikhalev, Shtern} over a field are Nielsen–Schreier. Recently, V. Dotsenko and U. Umirbaev \cite{Dots} showed that the varieties of all right-symmetric (or pre-Lie) algebras and all Lie-admissible algebras over a field of characteristic zero are also Nielsen–Schreier. As a consequence, every automorphism of a free left-symmetric algebra of finite rank over a field of characteristic zero is tame.

Let $A_n$ be a $n$-dimensional algebra with zero multiplication over a  field $K$ of characteristic $0$ with a linear basis $x_1,x_2,\ldots,x_n$. Its multiplicative universal enveloping algebra $U_n$ in the variety of left-symmetric algebras is a quadratic associative algebra with identity, generated by the operators of left multiplication $l_i=l_{x_i}$ and right multiplication $r_i=r_{x_i}$, where $1\leq i\leq n$. Moreover, $U_n$ contains both the polynomial algebra $L_n=K[l_1,\ldots,l_n]$ and the free associative algebra $R_n=K\langle r_1,\ldots,r_n\rangle$. 

This paper is devoted to the study of automorphisms and derivations of the algebra $U_n$. The main result of this paper establishes that the automorphism groups of the polynomial algebra $L_n$ and the algebra $U_n$ are isomorphic for all $n\geq 2$. Although we were not able to describe all locally nilpotent derivations, the proof of this isomorphism relies entirely on a detailed analysis of locally nilpotent derivations of a special form, which constitutes the most challenging part of our work. In contrast, we show that this isomorphism does not hold for $n=1$, and we provide a complete description of all automorphisms and locally nilpotent derivations of $U_1$.

The paper is organized as follows. In Section 2, we construct a linear basis and prove some preliminary results on the structure of $U_n$. In Section 3, we establish an embedding of   $\mathrm{Aut}(L_n)$ into $\mathrm{Aut}(U_n)$. Section 4 provides definitions and preliminary facts about derivations and locally nilpotent derivations. In Section 5, we prove that every locally nilpotent derivation of $U_n$ that induces the zero derivation on $L_n$ is itself zero for $n\geq 2$. In Section 6, we show that the embedding given in Section 3 is actually an isomorphism for all $n\geq 2$. Finally, in Section 7, we describe all locally nilpotent derivations and automorphisms of $U_1$.

\section{Structure of $U_n$}
\hspace*{\parindent}

Let $A$ be a left-symmetric algebra over an arbitrary field $K$. Denote by $U$ the universal (multiplicative) enveloping algebra 
\cite{Koz2} of  $A$ in the variety of left-symmetric algebras. Recall that $U$ is an associative unital algebra generated by the universal operators of left multiplication $l_x$ and right multiplication $r_x$, where $x\in A$. It is well known \cite{Jacobson} that the notion of an $A$-bimodule is equivalent to the notion of a left module over $U$. If $M$ is an $A$-bimodule then $M$ is a left $U$-module under the action 
\bes
l_x(a)=xa, r_x(a)=ax, \ \ \ a\in M, x\in A. 
\ees 
The identity \eqref{IN1} directly implies the defining relations 
of the algebra $U$:  
\begin{gather} \label{LT1}
l_xl_y-l_yl_x=l_{[x,y]},\;\; r_xl_y-l_yr_x-r_xr_y+r_{yx}=0,\;\; x,y\in A.
\end{gather}

D. Kozybaev and U. Umirbaev \cite{Koz2} constructed a linear basis of $U$.
\begin{theorem}\label{LTth1} \cite{Koz2}
Let $A$ be a left-symmetric algebra with a linear basis $x_1,x_2,\ldots,x_k,\ldots$. Then a linear basis of the universal 
enveloping algebra $U$ of $A$ consists of words of the form 
$$l_{x_{i_1}}l_{x_{i_2}}\ldots l_{x_{i_t}}r_{x_{j_1}}r_{x_{j_2}}\ldots r_{x_{j_s}},$$
where $i_1\leq i_2\leq \ldots \leq i_t$, $s,t\geq 0$.
\end{theorem}

Denote by $R$ the right universal enveloping algebra of $A$, i.e., the subalgebra of $U$ with identity $1$  
generated by the operators $r_x$, where $x\in A$. Similarly, we define the left universal enveloping algebra 
$L$  of $A$ as the subalgebra of $U$ with identity $1$ generated by the operators $l_x$, where $x\in A$.

\begin{corol} \label{LTcor1} \cite{Koz1} Under the conditions of Theorem \ref{LTth1}, the following statements are true: 

(a)  A linear basis of  $L$ consists of words of the form
$$l_{x_{j_1}}l_{x_{j_2}}\ldots l_{x_{j_t}}, \ \ j_1\leq j_2\leq \ldots \leq j_t, \ \ t\geq 0.$$

(b) The algebra $L$ is an associative unital algebra with generators $l_{x_i}$, $i\geq 1$, and defining relations
$$l_{x_j}l_{x_k}-l_{x_k}l_{x_j}-l_{[x_j,x_k]}=0, \;\; j>k.$$

(c) The algebra $R$ is a free associative algebra with free generators $r_{x_i}, \ i\geq 1$.
\end{corol}

Let $A_n$ be an $n$-dimensional algebra over  a field $K$ of characteristic zero, with zero multiplication  and a linear basis $x_1,x_2,\ldots,x_n$.
 The universal enveloping algebra $U_n$ of $A_n$ in the variety of left-symmetric algebras is generated by the operators $l_i=l_{x_i},r_i=r_{x_i}$, $1\leq i\leq n$.  The relations \eqref{LT1} imply the defining relations of $U_n$:
\begin{gather} \label{s1}
l_il_j=l_jl_i, \text{  } i>j,
\end{gather}
\begin{gather} \label{s2}
r_il_j=l_jr_i+r_ir_j.
\end{gather}

By Corollary \ref{LTcor1}, the left universal enveloping algebra
$L_n$ of $A_n$ is the polynomial algebra $L_n=K[l_1,\ldots, l_n]$ over $K$ in the variables $l_1,\ldots, l_n$, and 
the right universal enveloping algebra $R_n$ of $A_n$ is the free associative algebra $R_n=K\left\langle r_1,\ldots, r_n\right\rangle$ over $K$ with free generators $r_1,\ldots, r_n$. Consequently, the words of the form 
\bee\label{u1}
u=l_1^{s_1}l_2^{s_2}\ldots l_n^{s_n}, \ s_i\geq 0, 
\eee
form a linear basis of $L_n$ and the words of the form 
\bee\label{u2}
v=r_{j_1}r_{j_2}\ldots r_{j_s}
\eee
form a linear basis of $R_n$. 
By Theorem \ref{LTth1}, the basis of $U_n=L_nR_n$ consists of elements of the form 
\begin{gather} \label{LT3}
uv,
\end{gather}
where $u, v$ are words of the form (\ref{u1}) and (\ref{u2}), respectively.

Let $I_n$ be the ideal of $U_n$ generated by all $r_i$. Obviously, 
\bes
L_n\simeq U_n/I_n. 
\ees

Define the weight degree function $\mathrm{wdeg}$ on $U_n$ corresponding to the vector $w=(w_1,\\ w_2,\ldots,w_n)$, 
where $w_i$ are integers,  by setting $\mathrm{wdeg}(l_i)=$ $\mathrm{wdeg}(r_i)=w_i$, $1\leq i\leq n$. This degree function is compatible with relations \eqref{s1} and \eqref{s2}. Set $\mathrm{wdeg}(0)=-\infty$ as usual. 

Let us call an element $g\in U_n$ $w$-\emph{homogeneous} (\emph{homogeneous}) if all terms of $g$ have the same $w$-degree (degree) and let $(U_n)_m$ 
be the subspace of all $w$-homogeneous elements of the $w$-degree $m$. This determines the grading  
$$U_n=\oplus_{m\in \mathbb{Z}} (U_n)_m, \ \ \ \ (U_n)_s(U_n)_t\subseteq (U_n)_{s+t}.$$

The standard degree function $\deg$ corresponds to $w=(1,1,\ldots,1)$-weight degree function. We almost always use this degree function. The case $w=(1,\ldots,1,0)$ is used only once in the proof of Proposition \ref{th1}.

For any $s=(s_1,\ldots,s_n)\in \mathbb{Z}^n_+$, where $\mathbb{Z}_+$ is the set of non-negative integers, set $\left|s\right|=s_1+\ldots+ s_n$.    For any $s,t\in \mathbb{Z}^n_+$ set $s<t$ if $\left|s\right|<\left|  t\right|$ or $\left|s\right|=\left|t\right|$ and $s$ precedes $t$ with respect to the lexicographic order.

Let $u$ be a basis element of $L_n$ of form \eqref{u1}. Set $\mathrm{pdeg}(u)=(s_1,\ldots,s_n)$. If $u_1, u_2$ are two elements of the form \eqref{u1}, then set $u_1<u_2$ if  $\mathrm{pdeg}(u_1)< \mathrm{pdeg}(u_2)$. 

Arbitrary element $g\in U_n$ can be uniquely written as
$$g=u_1h_1+u_2h_2+\ldots+u_sh_s,$$ 
where $u_i$ are basis elements of $L_n$, $u_1>\ldots>u_n$, $h_i\in R_n$, $1\leq i\leq s$.
If $h_1\neq0$, put 
$$ \mathrm{Lm}(g)=u_1, \ \mathrm{Lc}(g)=h_1.$$  
Let us call $\mathrm{Lm}(g)$  and $\mathrm{Lc}(g)$ \emph{the leading monomial} and \emph{the leading coefficient} of $g$, respectively. For convenience, we set $\mathrm{Lm}(0)=\mathrm{Lc}(0)=0$.

Let $r_1>r_2>\ldots>r_n$. We define the deg-lex order on the set of basis elements of $R_n$ as in \cite [\S 8]{Bokut'}. That is, we say  $v_1<v_2$ for basis elements $v_1$ and $v_2$ of $R_n$  if one of the following conditions hold:
\begin{enumerate}
\item  $\deg(v_1)<\deg(v_2)$, 
\item $\deg(v_1)=\deg(v_2)$, $v_1$ precedes $v_2$ with respect to the lexicographic order.
\end{enumerate}

A linear map $D: U_n\rightarrow U_n$ is called \textit{ a derivation} of $U_n$, if it satisfies the Leibniz law, that is 
$$D(xy)=D(x)y+xD(y) \text{ for all } x,y\in U_n.$$ 

For any $a\in U_n$ the mapping $\mathrm{ad}_{a}:U_n\rightarrow U_n$ defined by $\mathrm{ad}_{a}(u)=[a,u]=au-ua$, where $u\in U_n$, is called an {\em inner} derivation. Using the inner derivations $\mathrm{ad}_{l_i}$, we can rewrite \eqref{s1} and \eqref{s2} as
\begin{gather} \label{form1}
\mathrm{ad}_{l_i}(l_j)=[l_i,l_j]=0, \ \ -\mathrm{ad}_{l_i}(r_j)=[r_j,l_i]= r_jl_i-l_ir_j=r_jr_i.
\end{gather}
Therefore, $\mathrm{ad}_{l_i}(v)\in R_n$ for all $v\in R_n$, and $\mathrm{ad}_{l_i}$ is a derivation of $R_n$ too. 

Denote by $\frac{\partial }{\partial l_j}$ the ordinary partial derivative of $L_n$ with respect to $l_j$.

\begin{lemma}\label{lem1}  Let $f(l_1,\dots,l_n)\in L_n$. Then the following statements are true:

$(a)$ $l_1-r_1, \dots, l_n-r_n$ generate the polynomial algebra $K[l_1-r_1, \dots, l_n-r_n]$;

$(b)$ $f(l_1,\dots,l_n)r_i=r_if(l_1-r_1,\dots,l_n-r_n)$;

$(c)$  $f(l_1-r_1,\dots,l_n-r_n)=f(l_1,\dots,l_n)-\sum^{n}_{j=1}\frac{\partial f(l_1,\dots,l_n)}{\partial l_j}r_j$.

\end{lemma}
\Proof 
 Write  \eqref{s2} as
\begin{gather} \label{s3}
l_jr_i=r_i(l_j-r_j). 
\end{gather}
Using this we have  $$ (l_i-r_i)(l_j-r_j)=l_il_j-l_ir_j-r_i(l_j-r_j)=l_jl_i-r_j(l_i-r_i)-l_jr_i=(l_j-r_j)(l_i-r_i).$$
Consequently, $l_1-r_1, \dots, l_n-r_n$  generate a commutative algebra. They generate the polynomial algebra $K[l_1-r_1, \dots, l_n-r_n]$ since their images in $U_n/I_n$ generate the polynomial algebra $L_n$. This proves $(a)$. 

Notice that $(b)$ is a direct consequence of $(a)$ and (\ref{s3}). 

It is sufficient to check $(c)$ for monomials $f\in L_n$. We prove this by induction on $\deg(f)=k$. Let $f=l_sg$.  
 By \eqref{s3} and by the induction proposition, 
$$f(l_1-r_1,\dots,l_n-r_n)=g(l_1-r_1,\dots,l_n-r_n) (l_s-r_s)=g(l_1,\dots,l_n)(l_s-r_s)$$
$$-\sum^{n}_{j=1}\frac{\partial g(l_1,\dots,l_n)}{\partial l_j}r_j(l_s-r_s)= g(l_1,\dots,l_n)l_s-g(l_1,\dots,l_n)r_s-\sum^{n}_{j=1}\frac{\partial g(l_1,\dots,l_n)}{\partial l_j}l_sr_j$$
$$=f(l_1,\dots,l_n)-g(l_1,\dots,l_n)\frac{\partial l_s}{\partial l_s}r_s-\sum^{n}_{j=1}\frac{\partial g(l_1,\dots,l_n)}{\partial l_j}l_sr_j$$
$$=f(l_1,\dots,l_n)-\sum^{n}_{j=1}\frac{\partial (l_sg(l_1,\dots,l_n))}{\partial l_j}r_j=f(l_1,\dots,l_n)-\sum^{n}_{j=1}\frac{\partial f(l_1,\dots,l_n)}{\partial l_j}r_j. \Box$$

\begin{corol}\label{c_1}  
Let $f(l_1,\dots,l_n)\in L_n$. Then  $$r_if(l_1,\dots,l_n)=f(l_1,\dots,l_n)r_i+r_i \sum^{n}_{j=1}\frac{\partial f(l_1,\dots,l_n)}{\partial l_j}r_j. $$
\end{corol}
\Proof By Lemma \ref{lem1}, we have
$$r_if(l_1,\dots,l_n)=r_if(l_1-r_1,\dots,l_n-r_n)+r_i\sum^{n}_{j=1}\frac{\partial f(l_1,\dots,l_n)}{\partial l_j}r_j$$
$$=f(l_1,\dots,l_n)r_i+r_i\sum^{n}_{j=1}\frac{\partial f(l_1,\dots,l_n)}{\partial l_j}r_j. \Box$$

\begin{lemma}\label{LTlem3} 
Let $w_1=u_1v_1$ and
$w_2=u_2v_2$ be basis elements of $U_n$ of the form (\ref{LT3}). Then 
$$w_1w_2=u_1u_2v_1v_2+ h,$$
where $\mathrm{Lm} (h)<u_1u_2$.
\end{lemma}
\Proof It suffices to consider the case $w_1=r_{i_1}\ldots r_{i_p}$ and $w_2=l_j$. We have
$$w_1w_2= r_{i_1} \ldots r_{i_p}l_j=l_j r_{i_1} \ldots r_{i_p}-[l_j,r_{i_1} \ldots r_{i_p}]
= l_j r_{i_1} \ldots r_{i_p}-ad_{l_j}(r_{i_1} \ldots r_{i_p}).$$
By \eqref{form1}, $\mathrm{ad}_{l_j}(r_{i_1} \ldots r_{i_p})\in R_n$. $\Box$

\begin{corol}\label{LTcor2}  The algebra $U_n$ has no zero divisors and 
\bes
\mathrm{Lm}(g\cdot h)=\mathrm{Lm}(g)Lm(h), 
\ees
\bes
\mathrm{wdeg}(g\cdot h)=\mathrm{wdeg} (g)+\mathrm{wdeg}(h), 
\ees
for all $g,h\in U_n$. 
\end{corol}
\Proof Lemma \ref{LTlem3} immediately implies that $U_n$ 
  has no zero divisors and establishes the first equality of the corollary. This, in turn, implies the second equality as well. $\Box$

\begin{lemma}\label{lem_1}   The element $ r_i^kr_jh$, where $k\geq 1$, $i\neq j$, $h\in R_n$, can be written as
$$r_i^kr_jh=\mathrm{ad}_{l_i}(r_iu)+r_ir_jv$$
for some $u,v\in R_n$.
\end{lemma}
\Proof We prove the statement of the lemma by induction on $k$. Let $k=1$. Then $r_ir_jh=\mathrm{ad}_{l_i}(r_iu_0)+r_ir_jv_0$, where $u_0=0$ and $v_0=h$.
If $k>1$ then 
$$r_i^kr_jh=\frac{1}{k-1}\mathrm{ad}_{l_i}(r_i^{k-1})r_jh=\frac{1}{k-1}\mathrm{ad}_{l_i}(r_i^{k-1}r_jh)-\frac{1}{k-1}r_i^{k-1}\mathrm{ad}_{l_i}(r_jh).$$
Notice that $\mathrm{ad}_{l_i}(r_jh)=r_jh_0$ for some $h_0\in R_n$. Applying the induction proposition to $r_i^{k-1}\mathrm{ad}_{l_i}(r_jh)=r_i^{k-1}r_jh_0$, we obtain 

$$r_i^kr_jh=\frac{1}{k-1}\mathrm{ad}_{l_i}(r_i^{k-1}r_jh)+\mathrm{ad}_{l_i}(r_iu_{k-1})+r_ir_jv_{k-1}=\mathrm{ad}_{l_i}(r_iu_k)+r_ir_jv_k. \Box$$

\begin{lemma}\label{NLNlem2} 
Let $u_1,u_2, \ldots, u_n \in I_n$, $n\geq 2$, such that
\begin{equation}\label{eq_11}
\mathrm{ad}_{l_j}(u_i)=\mathrm{ad}_{l_i}(u_j). 
\end{equation} 
Then $u_i=\mathrm{ad}_{l_i} (g)$ for some  $ g\in I_n.$ 
\end{lemma}

\Proof Since the equations \eqref{eq_11} are homogeneous, we can assume that all $u_i$, $1\leq i\leq n$, are homogeneous of the same degree $t$. If $t=1$, then $u_i\in R_n$ and $\mathrm{ad}_{l_j}(u_i)=u_ir_j$ for all $i,j$. Consequently, \eqref{eq_11} holds only in the case $u_i=0$ for all $i$, which contradicts to $t=1$. 

Suppose $t>1$. 
First consider the case  $u_i\in R_n$ for all $1\leq i \leq n$. Let $u_i=r_1u_{i1}+r_2u_{i2}+\ldots+r_nu_{in}$.  Assume that $u_{ii}=\sum_{i=1}^nr_kc_k$. 

Notice that if $m=r_i^t$ then $r_i^2m=-\frac{1}{t+1}\mathrm{ad}_{l_i}(r_i^{t+1})$ and if $m=r_i^tr_jm_1$, $i\neq j$, then $r_i^2m=\mathrm{ad}_{l_i}(r_iu)+r_ir_jv$ for some $u,v\in R_n$ by Lemma \ref{lem_1}. Then 
\bes
r_i^2c_i=\sum^n_{k=1,k\neq i}r_ir_ka_{ik}'+ad_{l_i}(r_ib_i), 
\ees
for some $a_{ij}',b_i\in R_n$. Consequently, for any $i$ we obtain 
\begin{equation}\label{equ_2}
 r_iu_{ii}=\sum^n_{k=1,k\neq i}r_ir_ka_{ik}+ad_{l_i}(r_ib_i),
\end{equation}
for some $a_{ij},b_i\in R_n$. 

For $j=1, i\neq 1$ the equation \eqref{eq_11} implies 
$$\mathrm{ad}_{l_1}(r_1u_{i1}+\ldots+r_{i-1}u_{i(i-1)}+ \sum^n_{k=1, k\neq i}r_ir_ka_{ik}+\mathrm{ad}_{l_i}(r_ib_i)+r_{i+1}u_{i(i+1)}+ \ldots + r_nu_{in} )$$
$$=\mathrm{ad}_{l_i}(\sum^n_{k=2}r_1r_ka_{1k}+\mathrm{ad}_{l_1}(r_1b_1)+r_2u_{12}+r_3u_{13}+ \ldots + r_nu_{1n}).$$
Comparing the associative words that start to $r_1$ and $r_i$ we obtain the following two equations: 
$$\mathrm{ad}_{l_1}(r_1u_{i1}-\mathrm{ad}_{l_i}(r_1b_1)) =\mathrm{ad}_{l_i}(\sum^n_{k=2}r_1r_ka_{1k}),$$
$$\mathrm{ad}_{l_i}(r_iu_{1i}-\mathrm{ad}_{l_1}(r_ib_i))=\mathrm{ad}_{l_1}(\sum^n_{k=1, k\neq i}r_ir_ka_{ik}).$$
 If $r_1u_{i1}-\mathrm{ad}_{l_i}(r_1b_1)=r_1w_1\neq 0$ and $r_iu_{1i}-\mathrm{ad}_{l_1}(r_ib_i)=r_iw_2\neq 0$, then on the left sides of the last two equalities there are  nonzero elements of the form $r_1r_1h$ and $r_ir_iq$, respectively. But there are no such elements on the right sides of the last two equalities. Therefore $r_1u_{i1}=\mathrm{ad}_{l_i}(r_1b_1)$ and  $\sum^n_{k=1, k\neq i}r_ir_ka_{ik}=0$. Consequently, by \eqref{equ_2}, $$ r_iu_{ii}=\mathrm{ad}_{l_i}(r_ib_i).$$
Similarly, we prove $$r_ju_{ij}=\mathrm{ad}_{l_i}(r_jb_j) \text{ for any }  j\neq i.  $$
Therefore, 
 $$u_i=r_1u_{i1}+r_2u_{i2}+\ldots+r_nu_{in}=\mathrm{ad}_{l_i}(r_1b_1+r_2b_2+\ldots+r_nb_n) \ \text{ for all } 1\leq i \leq n.$$

Let
$$u_i=f_{i1}v_{i1}+\ldots+ f_{in}v_{in}\in I_n,$$
where $f_{ij}$ are basis elements of $L_n$, $f_{i1}>\ldots >f_{in}$, and $v_{ij}\in R_n\cap I_n$ for all $1\leq i,j \leq n$.
By \eqref{eq_11},
$$ f_{i1}\mathrm{ad}_{l_j}(v_{i1})+\ldots+ f_{in}\mathrm{ad}_{l_j}(v_{in})=f_{j1}\mathrm{ad}_{l_i}(v_{j1})+\ldots+ f_{jn}\mathrm{ad}_{l_j}(v_{jn}).$$
It follows that 
$f_{ik}=f_{jk}$ and  $\mathrm{ad}_{l_j}(v_{ik})=\mathrm{ad}_{l_i}(v_{jk})$ for all $1\leq k \leq n.$ As proven above, there exist $g_k\in R_n$ such that $v_{ik}=\mathrm{ad}_{l_i}(g_k)$ for all $1\leq i,k \leq n$. Consequently
$$ u_i=f_{1}v_{i1}+\ldots+ f_{n}v_{in}=\mathrm{ad}_{l_i}(f_1g_1+\ldots+f_ng_n). \Box$$

\begin{lemma}\label{U_14} Let $g\in I_n$. If 
$-\mathrm{ad}_{l_i}(g)=r_i g+gr_i$,  
 then  $\mathrm{Lc}(g)=\alpha_1 r_ir_1+\alpha_2 r_i r_2+\ldots+\alpha_n r_i r_n$.
\end{lemma}
\Proof Let $\mathrm{Lm}(g)=f$ and $\mathrm{Lc}(g)=u$. Then $g=fu+w$,  $w\in I_n$  and $\mathrm{Lm}(w)<f$. By the condition of lemma, we have 
$$ -f\mathrm{ad}_{l_i}(u)-\mathrm{ad}_{l_i}(w)=r_ifu+r_iw+fur_i+wr_i.$$
 By Lemma \ref{LTlem3}, this implies $$-\mathrm{ad}_{l_i}(u)=r_iu+ur_i.$$

 Since $u\in R_n$ we can write $u=u_1r_1+\ldots+u_nr_n$, where $u_j\in R_n$ for all $j$. Then the last equality can be rewritten as 
$$ -\mathrm{ad}_{l_i}(u_1)r_1+\ldots-\mathrm{ad}_{l_i}(u_n)r_n+u_1r_1r_i+\ldots+u_nr_nr_i=r_iu_1r_1+\ldots+r_iu_nr_n+u_1r_1r_i+\ldots+u_nr_nr_i,$$
and consequently, 
$$ -\mathrm{ad}_{l_i}(u_1)r_1+\ldots-\mathrm{ad}_{l_i}(u_n)r_n=r_iu_1r_1+\ldots+r_iu_nr_n.$$
Therefore 
 $$ -\mathrm{ad}_{l_i}(u_j)=r_iu_j  \ \text{ for all } \  1\leq j\leq n.$$
It is possible only if $u_j=r_iv_j$  for some $v_j$ and the last equation takes the form $r_ir_iv_j-r_i\mathrm{ad}_{l_i}(v_j)=r_ir_iv_j$ for all $j$. Hence $\mathrm{ad}_{l_i}(v_j)=0$. Consequently, $v_j\in K$,  $u_j=\alpha_j r_i$ and 
$$u=\alpha_1 r_ir_1+\ldots+\alpha_n r_ir_n,$$
where  $\alpha_j \in K$ for all  $j$. $\Box$

\section{Polynomial automorphisms}
\hspace*{\parindent}

 Notice that 
\bes 
U_n=L_n\oplus I_n, 
\ees 
is a direct sum of vector spaces. Then every element $q\in U_n$  can be uniquely written as
$$q=q(l)+q_0,$$
where $q(l)\in L_n$, $q_0\in I_n$.

Denote by $\varphi=(g_1,\ldots,g_n,h_1,\ldots,h_n)$ the  endomorphism of $U_n$, such that $\varphi(l_i)=g_i$ and $\varphi(r_i)=h_i$ for all $1\leq i\leq n$. 

\begin{lemma}\label{ELlem1} 
Let $\varphi=(g_1,\ldots,g_n,h_1,\ldots,h_n)$ be an arbitrary endomorphism of $U_n$. Then $\varphi(I_n)\subseteq I_n$.
\end{lemma}
\Proof Applying the endomorphism $\varphi$ to the relation \eqref{s2} in the case $i=j$, we obtain
$$(h_i(l)+h_{i0})(g_i(l)+g_{i0})=(g_i(l)+g_{i0})(h_i(l)+h_{i0})+(h_i(l)+h_{i0})(h_i(l)+h_{i0}).$$
From this, we can deduce that $h_i(l)g_i(l)=g_i(l)h_i(l)+h_i(l)h_i(l)$ and $h_i(l)^2=0$.  Hence, $h_i(l)=0$. $\Box$

Denote by $\mathrm{Aut}(L_n)$ and $\mathrm{Aut}(U_n)$  the groups of all automorphisms of $L_n$ and $U_n$, respectively. For any 
$n$-tuple $\varphi=(f_1,\ldots, f_n)$ of elements of $L_n$ set 
\begin{equation}\label{equ1} 
\varphi^*=(f_1,\ldots, f_n, \sum^{n}_{s=1}\frac{\partial f_1}{\partial l_s}r_s,\ldots, \sum^{n}_{s=1}\frac{\partial f_n}{\partial l_s}r_s). 
\end{equation}

\begin{pr}\label{up1} If $\varphi$ is an automorphism of $L_n$ then $\varphi^*$ is an automorphism of $U_n$, and the map 
\bes
\Phi : \mathrm{Aut}(L_n)\rightarrow \mathrm{Aut}(U_n)
\ees
defined by $\Phi(\varphi)=\varphi^*$ is an embedding of groups.
\end{pr} 
\Proof First of all, let's show that $\varphi^*$ is an endomorphism of $U_n$. Indeed, applying $\varphi^*$ to the relation \eqref{s2}, we obtain
\bes
(\sum^{n}_{s=1}\frac{\partial f_i}{\partial l_s}r_s)f_j
=f_j(\sum^{n}_{s=1}\frac{\partial f_i}{\partial l_s}r_s) +(\sum^{n}_{s=1}\frac{\partial f_i}{\partial l_s}r_s)(\sum^{n}_{s=1}\frac{\partial f_j}{\partial l_s}r_s).
\ees
By Corollary \ref{c_1}, this equality holds.

 Let $\phi=(f_1,\ldots, f_n), \ \psi=(g_1,\ldots, g_n)$ be any automorphisms of $L_n$. Denote by $g_i(f)$ the element $g_i(f_1, \ldots, f_n)$. We have
$$\Phi(\phi\psi)= \Phi ((g_1(f),\ldots, g_n(f)))$$
$$= (g_1(f),\ldots, g_n(f),\sum^{n}_{s=1}\frac{\partial g_1(f)}{\partial l_s}\sum^{n}_{t=1}\frac{\partial f_s}{\partial l_t}r_t, \ldots,\sum^{n}_{s=1}\frac{\partial g_n(f)}{\partial l_s}\sum^{n}_{t=1}\frac{\partial f_s}{\partial l_t}r_t) $$
$$=(f_1,\ldots, f_n, \sum^{n}_{t=1}\frac{\partial f_1}{\partial l_t}r_t,\ldots, \sum^{n}_{t=1}\frac{\partial f_n}{\partial l_t}r_t)(g_1,\ldots, g_n, \sum^{n}_{s=1}\frac{\partial g_1}{\partial l_s}r_s,\ldots, \sum^{n}_{s=1}\frac{\partial g_n}{\partial l_s}r_s)$$
$$ =\Phi(\phi)\Phi(\psi).$$
It is clear that $\Phi(id)=id$, where $id$ is the identity endomorphism. Consequently, $\Phi$ sends automorphisms to automorphisms and is an  embedding of groups. $\Box$

Let $\mathrm{Aff}(L_n)$ be the group of all affine automorphisms of $L_n$. Denote by $\mathrm{Aff}(U_n)$ the group of all  automorphisms $\psi$ of   $U_n$ such that $\deg(\psi(l_i))=\deg(\psi(r_i))=1$ for all $1\leq i\leq n$.

\begin{lemma} \label{affaut} $\Phi(\mathrm{Aff}(L_n))=\mathrm{Aff}(U_n)$ if $n\geq 2$. 
\end{lemma}
\Proof  Let $\theta\in \mathrm{Aff}(U_n)$. By Lemma \ref{ELlem1},  $\theta$ induces an automorphism $\phi$ of $L_n$. Obviously, 
$\phi^*\in \mathrm{Aff}(U_n)$. Changing $\theta$ by $\theta(\phi^*)^{-1}$, we may assume that 
  $\theta(l_i)=l_i+w_i$ and $\theta(r_i)=h_i $, where   $w_i,h_i$ are linear elements of $R_n$.  Applying $\theta$ to the relation \eqref{s2}, we obtain 
$ h_il_j+h_iw_j=l_jh_i+w_jh_i+h_ih_j.$
 By \eqref{s2}, we can write this as $$l_jh_i+h_ir_j+h_iw_j=l_jh_i+w_jh_i+h_ih_j$$ since $h_i, w_j$ are linear elements of $R_n$.
Hence
$$h_i(r_j+w_j-h_j)=w_jh_i.$$
It follows that $w_j=\beta_{ji}h_i$, $\beta_{ji}\in K$, and $h_j=r_j$. Therefore  $w_j=\beta_{ji}r_i$  for all $i,j$.  Consequently, $w_j=0$ for all $j$. $\Box$

\section{Derivations of $U_n$}
\hspace*{\parindent}

Let $\mathrm{Der}\, U_n$ be the Lie algebra of all derivations of $U_n$.  Denote by 
\begin{gather} \label{GR1}
D=\sum^{n}_{i=1}u_i \frac{\partial}{\partial l_i}+\sum^{n}_{i=1}v_i \frac{\partial}{\partial r_i}
\end{gather}
the derivation of $U_n$ such that $D(l_i)=u_i$, $D(r_i)=v_i$ for all $1\leq i \leq n$. 

\begin{lemma}\label{TRlem1} Let $D$ be a derivation of $U_n$. Then $D(I_n)\subseteq I_n$.
\end{lemma}
\Proof Let $D$ has the form (\ref{GR1}). 
 Applying $D$ to the relation 
\eqref{s2} when $i=j$, we obtain 
\begin{gather} \label{TD0}
v_il_i+r_iu_i=u_ir_i+l_iv_i+v_ir_i+r_iv_i.
\end{gather} 

Denote by $J$ the ideal of $U_n$ generated by all $r_ir_j$ for all  $1\leq i,j \leq n$. Let $g=g(l)+g_0\in U_n$, where  $q(l)\in L_n$, $q_0\in I_n$. By Lemma \ref{LTlem3} and Corollary \ref{c_1}, we have 
$$gl_i=g(l)l_i+g_0l_i=g(l)l_i+l_ig_0+q_1,  \ q_1\in J,$$
 $$r_ig=r_ig(l)+r_ig_0=g(l)r_i+q_2, \ q_2\in J.$$
Using this and \eqref{TD0}, we obtain
$$v_i(l)l_i+l_iv_{i0}+u_i(l)r_i+q_3=u_i(l)r_i+l_iv_i(l)+l_iv_{i0}+v_i(l)r_i+v_i(l)r_i+q_4,  $$
where $q_3,q_4 \in J$. Consequently, $2v_i(l)r_i+q=0$,   where $q\in J$. Hence, $ v_i(l)=0$. $\Box$

Consequently, every derivation of $U_n$ induces a derivation of $L_n$. But we do not know if there is a way to continue any derivation of $L_n$ to a derivation of $U_n$. Moreover, there is an example of a nonzero derivation of $U_n$ that induces a zero derivation of $L_n$. 

\begin{ex} \label{ex1} There is a derivation $D$ of $U_2$ such that
$$D(l_1)=r_1 r_1, \ \ D(l_2)=r_1 r_2, \ \ D(r_1)=0, \ \ D(r_2)=r_1 r_2-r_2 r_1.$$
\end{ex} 
In fact, direct calculations give that 
$$D(l_1l_2-l_2l_1)=r_1r_1l_2+l_1r_1r_2-r_1r_2l_1-l_2r_1r_1 $$
$$=l_2r_1r_1+r_1r_2r_1+r_1r_1r_2+l_1r_1r_2-l_1r_1r_2-r_1r_1r_2-r_1r_2r_1-l_2r_1r_1=0,$$
$$D(r_1l_1-l_1r_1-r_1r_1)= r_1r_1r_1- r_1r_1r_1=0,$$
$$D(r_1l_2-l_2r_1-r_1r_2)= r_1r_1r_2-r_1r_2r_1-r_1(r_1 r_2-r_2 r_1)=0,$$
$$ D(r_2l_1-l_1r_2-r_2r_1)=(r_1 r_2-r_2 r_1)l_1+r_2r_1 r_1-r_1 r_1r_2-l_1(r_1 r_2-r_2r_1)-(r_1 r_2-r_2 r_1)r_1$$
$$ =r_1r_1r_2-r_2r_1r_1+r_1r_2r_1-r_2r_1r_1+2r_2r_1 r_1- r_1 r_1r_2-r_1 r_2r_1=0,$$
$$D(r_2l_2-l_2r_2-r_2r_2)= (r_1 r_2-r_2 r_1)l_2+r_2r_1 r_2-r_1 r_2r_2-l_2(r_1 r_2-r_2r_1)-(r_1 r_2-r_2 r_1)r_2$$
$$ -r_2(r_1 r_2-r_2 r_1)=r_1 r_2r_2-r_2 r_2r_1+r_1r_2r_2-r_2r_1r_2-2r_1 r_2r_2+r_2 r_1r_2+r_2r_2 r_1=0.$$
Consequently, $D$ preserves all defining relations of $U_2$.

Recall that a derivation $D$ of $U_n$ is called \textit{locally nilpotent} if for every $a\in U_n$ there exists a natural number $m=m(a)$  such that $D^m(a)=0$. 
Notice that the derivation $D$ from Example \ref{ex1} is not locally nilpotent since $D^k(r_2)=\mathrm{ad}(r_1)^k(r_2)\neq 0$. 

Let $w=(w_1,\ldots,w_n)\in \mathbb{Z}^n$. 
Recall that $(U_n)_m$ denotes  be the subspace of all $w$-homogeneous elements of $w$-degree $m$ of $U_n$. Every nonzero element of $U_n$ can be written as 
$$g=g_r+g_{r+1}+\cdots+g_s\in U_n,$$
 where $g_i\in (U_n)_i$, $r\leq i\leq s$ and $g_s\neq 0$. In this case we call $\widehat{g}=g_s$ the {\em highest homogeneoous part} of $g$.

Notice that a formal expression of the form (\ref{GR1}) is not necessarily a derivation of $U_n$. For any formal expression of the forms 
\bes
u\frac{\partial}{\partial l_i}, u\frac{\partial}{\partial r_i},  u\in (U_n)_m, 
\ees
set 
\bes
\mathrm{wdeg}(u\frac{\partial}{\partial l_i})=\mathrm{wdeg}(u\frac{\partial}{\partial r_i})= m-w_i. 
\ees
We say that a formal expression $D$ of the form (\ref{GR1}) is {\em $w$-homogeneous of $w$-degree $m$} if all nonzero expressions $u_i \frac{\partial}{\partial l_i},v_i \frac{\partial}{\partial r_i}$ are $w$-homogeneous of $w$-degree $m$.

Fortunately, for any derivation $D$ of $U_n$, all of its homogeneous components are themselves derivations of $U_n$, since all defining relations of $U_n$ are $w$-homogeneous. 
Let $\mathrm{Der}_m U_n$ be the subspace of all $w$-homogeneous elements  of $\mathrm{Der}\,U_n$ of $w$-degree $m$.  This gives a grading 
$$\mathrm{Der}\,U_n=\oplus_{m\in  \mathbb{Z}} \mathrm{Der}_m U_n$$
of the Lie algebra of derivations.  
Every derivation can be written as 
$$D=D_p+D_{p+1}+\ldots+D_q,$$
where $D_i\in Der_iU_n$, $p\leq i\leq q$, and $D_q\neq 0$. Then $\widehat{D}=D_q$  is called the {\em highest homogeneous part} of $D$. 

The statements of the next two lemmas are well known (see, for example,\cite{M-L2}).

\begin{lemma}\label{GRlem1} 
$(\mathrm{Der}_m U_n)((U_n)_k)\subseteq (U_n)_{m+k}$.
\end{lemma}
\Proof Let $D\in \mathrm{Der}_m U_n$ be a derivation of the form \eqref{GR1}.
Then $u_i, v_i \in (U_n)_{m+w_i}$ for all $1\leq i\leq n$. 

Let $a\in (U_n)_k$ be a basis element of $U_n$ of the form (\ref{LT3}). By induction on $\deg(a)$ we prove that $D(a)\in (U_n)_{m+k}$. If $\deg a=1$, then $a=l_j$ or $a=r_j$ and $D(a)=u_j$ or $D(a)=v_j$, respectively. 
 If $\deg a>1$, then $a=a_1a_2$, where $\deg a_1=s\geq 1$, $\deg a_2=t\geq 1$. By the induction proposition, $D(a_1)\in (U_n)_{m+s}$ and $D(a_2)\in (U_n)_{m+t}$. Consequently, $D(a)=D(a_1a_2)=D(a_1)a_2+a_1D(a_2)\in (U_n)_{m+s+t}=(U_n)_{m+k}$. 
 $\Box$

\begin{lemma}\label{GRpr2} 
The highest homogeneous part of any locally nilpotent derivation of $U_n$ is locally nilpotent. 
\end{lemma}
\Proof Using Lemma \ref{GRlem1}, it is easy to check that 
\bes
\widehat{D(g)}=\widehat{D}(\widehat{g}) 
\ees
for any derivation $D$ and $g\in U_n$ if $\widehat{D}(\widehat{g})\neq 0$. This implies that $\widehat{D}$ is locally nilpotent whenever $D$ is. $\Box$

If $D$ is any formal expression of the form (\ref{GR1}) and $f\in L_n$, then we set 
\bes
fD=\sum^{n}_{i=1}fu_i \frac{\partial}{\partial l_i}+\sum^{n}_{i=1}fv_i \frac{\partial}{\partial r_i}. 
\ees
This defines an action of $L_n$ on the set of all formal expressions of the form (\ref{GR1}). 

For any formal expression $D$ of the form (\ref{GR1}) set 
\bes
D^{(r)}=\sum^{n}_{i=1}v_i \frac{\partial}{\partial r_i}. 
\ees

A formal expression $D$ of the form (\ref{GR1}) will be called {\em purely associative} if $D(l_i)=u_i\in R_n$, $D(r_i)=v_i\in R_n$ for all $1\leq i \leq n$.

Every nonzero derivation $D\in Der\,U_n$ can be uniquely written as
\bee\label{u3}
D=g_1D_1+g_2D_2+\ldots+g_mD_m,
\eee
where $g_i$ are monomials of $L_n$ and $D_i$ are nonzero purely associative formal expressions of the form \eqref{GR1}  for all $i$,  and $g_1>\ldots>g_m$.  Set also 
$$\mathrm{Lm}(D)=g_1, \ \mathrm{Lc}(D)=D_1.$$
Let us call $\mathrm{Lm}(D)$ and $\mathrm{Lc}(D)$ \emph{the leading monomial} and \emph{the leading coefficient} of $D$. For convenience, we set $\mathrm{Lm}(0)=\mathrm{Lc}(0)=0$.

\begin{lemma} \label{NLNpr2} 
 If $D$ is a locally nilpotent derivation of $U_n$ then $\mathrm{Lc}(D)^{(r)}$ is a locally nilpotent derivation of $R_n$.
\end{lemma}
\Proof 
Let $D$ be a locally nilpotent derivation of $U_n$ of the form (\ref{u3}). Then $\mathrm{Lm}(D)=g_1=g$, $\mathrm{Lc}(D)=D_1$, and $\mathrm{Lc}(D)^{(r)}=D_1^{(r)}$. We check the equality 
\bee \label{equ3}
\mathrm{Lc}(D(b))=D_1^{(r)}(\mathrm{Lc}(b))
\eee
for any $b\in U_n$ if $D_1^{(r)}(\mathrm{Lc}(b))\neq 0$. We may assume that $b$ is a basis element of $U_n$, i.e., $b=ha$, where $h$ is monomial of $L_n$, $a=r_{i_{1}}\ldots r_{i_{k}}$. Then \eqref{equ3} can be expressed in the form 
\bes
D(b)=gh \cdot (D_1^{(r)}(\mathrm{Lc}(b)))+w, \ \ \mathrm{Lm}(w)<gh.
\ees

By Lemma \ref{LTlem3}, we have
$$h D(a)=hD(r_{i_{1}}\ldots r_{i_{k}})=h(D(r_{i_{1}})r_{i_{2}}\ldots r_{i_{k}}+ \ldots
+r_{i_{1}}r_{i_{2}}\ldots D(r_{i_{k}}))$$
$$=hg(D_1(r_{i_{1}})r_{i_{2}}\ldots r_{i_{k}}+\ldots
+r_{i_{1}}r_{i_{2}}\ldots D_1(r_{i_{1}}))+w=hg (D_1^{(r)}(a))+w,$$
where $\mathrm{Lm}(w)<hg$. It is clear that 
\bes
\deg(\mathrm{Lm}(h D(a)))=\deg(h)+\deg(g)>\deg(\mathrm{Lm}( D(h)a))\geq \deg(h)+\deg(g)-1  
\ees
since $g\geq \mathrm{Lm}(D(l_i))$ for all $i$. Consequently, 
$$D(b)=gh\cdot(D_1^{(r)}(a))+w=gh\cdot(D_1^{(r)}(\mathrm{Lc}(b)))+w,$$
where $\mathrm{Lm}(w)<gh$.

By \eqref{equ3}, if $c\in R_n$ and $(D_1^{(r)})^k(c)\neq 0$, then 
\bes
\mathrm{Lc}(D^k(c))=(D_1^{(r)})^k(c), 
\ees
for all $k\geq 0$. Therefore, $D_1^{(r)}$ is locally nilpotent because $D$ is. $\Box$

\section{Locally nilpotent derivations of $U_n$, $n\geq 2$}
\hspace*{\parindent}

The main result of this section, which will serve as our primary tool for studying automorphisms, is the following theorem.

\begin{theorem} \label{tu1}
Every locally nilpotent derivation of $U_n$ that induces the zero derivation on $L_n$ is itself the zero derivation. 
\end{theorem}
\Proof If $D$ induces the zero derivation on $L_n$ then $D(U_n)\subseteq I_n$. 
Suppose that $D$ is a nonzero derivation of $U_n$ satisfying the conditions of the theorem. By Lemma \ref{GRpr2}, its highest homogeneous part also satisfies the conditions of the theorem. So we may assume that $D$ is homogeneous with respect to the standard degree function $\deg$.

Suppose that $D$ is written in the form (\ref{u3}),  $\mathrm{Lm}(D)=g_1=g$, 
and 
\bes
\mathrm{Lc}(D)=D_1=\sum^{n}_{i=1}u_i \frac{\partial}{\partial l_i}+\sum^{n}_{i=1}v_i \frac{\partial}{\partial r_i}, 
\ees
where $u_i, v_i\in R_n\cap I_n$ for all $1\leq i \leq n$. This means that 
\begin{equation} \label{NLN6}
\begin{split}
D(l_i)=g u_{i}+a_{i}, \ \ D(r_i)=g v_{i}+b_{i}, 
\end{split}
\end{equation}
where  $a_i,b_i\in I_n$, $\mathrm{Lm} (a_i), \mathrm{Lm}(b_i)<g$ for all $i$. By Lemma \ref{NLNpr2}, $D_1^{(r)}$ is a locally nilpotent derivation of the associative algebra $R_n$. 

Applying $D$ to the relations  \eqref{s2}, we obtain 
$$ (g v_{i}+b_{i})l_j+r_i(g u_{j}+a_{j})=(g u_{j}+a_{j})r_i+l_j(g v_{i}+b_{i})+(g v_{i}+b_{i})r_j+r_i(g v_{j}+b_{j})$$
for all $i,j$. 
By Lemma \ref{LTlem3}, this can be rewritten as 
$$ g  l_j v_{i} -g \mathrm{ad}_{l_j}(v_{i})+l_j b_{i}-ad_{l_j}(b_i)+g r_i u_{j}+c_1$$
$$= g u_{j}r_i+g l_jv_{i}+l_jb_{i}+g v_{i}r_j+b_{i}r_j+g r_iv_{j}+c_2,$$
where $\mathrm{Lm}(c_1), \mathrm{Lm}(c_2)<g$. Consequently, 
\bes
g(-\mathrm{ad}_{l_j}(v_{i})-v_ir_j)+g(r_i u_{j}-u_{j}r_i)+c
= g r_iv_{j},
\ees
where $\mathrm{Lm}(c)<g$. Comparing the leading monomials, from this we obtain 
\begin{gather} \label{NLN7}
 -\mathrm{ad}_{l_j}(v_{i})-v_ir_j+r_i u_{j}-u_{j}r_i= r_iv_{j}. 
\end{gather}

Applying $D$ to the relations \eqref{s1}, we obtain $D(l_i)l_j+l_iD(l_j)=D(l_j)l_i+l_jD(l_i)$ for all $i<j$. This implies $\mathrm{ad}_{l_i}(D(l_j))=\mathrm{ad}_{l_j}(D(l_i)).$ By Lemma \ref{NLNlem2}, there exists  $w\in I_n$ such that $D(l_i)=\mathrm{ad}_{l_i} (w)$ for all $i$. Consequently, $\mathrm{Lm}(D(l_i))=f=\mathrm{Lm}(w)$ for all $1 \leq i \leq n$.

Obviously, $g\geq f$. The remaining part of the proof we separate into the following 3 cases: 

{\em Case 1}. $g>f$. 

{\em Case 2}. $g=f$ and $v_1=0$. 

{\em Case 3}. $g=f$ and $v_i\neq 0$ for all $i$. 

All cases are shown to be impossible in Lemmas \ref{NLNlem8}, \ref{NLNlem10}, and \ref{NLNlem9}, respectively. 
Specifically, each lemma demonstrates $D_1^{(r)}$ cannot be locally nilpotent. 
This contradiction implies $D=0$. $\Box$

The following three lemmas are formulated under the assumptions established in the proof of Theorem \ref{tu1}.

\begin{lemma}\label{NLNlem8} If $g>f$, then $D_1^{(r)}$ is not locally nilpotent. 
\end{lemma} 
\Proof The condition $g>f$ means that $u_1=\ldots=u_n=0$. Then the equality \eqref{NLN7} implies 
\begin{gather} \label{eq_16}
-\mathrm{ad}_{l_j}(v_i)-v_ir_j=r_i v_{j}.
\end{gather}
 If $j=i$ in \eqref{eq_16}, by Lemma \ref{U_14}, we obtain
  $$v_{i}=\alpha_{i1} r_i r_1+\alpha_{i2} r_i r_2+\ldots+\alpha_{in} r_i r_n$$
for all $1 \leq i \leq n$, where $\a_{ij}\in K$.
Substituting $v_i,v_j$, $i\neq j$, into \eqref{eq_16}, we obtain 
$$\alpha_{i1} r_ir_j r_1+\alpha_{i2} r_i r_jr_2+\ldots+\alpha_{in} r_i r_jr_n=\alpha_{j1}r_i r_j r_1+\alpha_{j2} r_ir_j r_2+\ldots+\alpha_{jn} r_ir_j r_n.$$
Hence $\alpha_{ik}=\alpha_{jk}$ for all $1\leq k\leq n$. Thus 
$$v_{i}=\alpha_{1} r_i r_1+\alpha_{2} r_i r_2+\ldots+\alpha_{n} r_i r_n, \ 1 \leq i \leq n,$$ 
and 
 $$D_1^{(r)}=\sum^{n}_{i=1}r_i\left (\alpha_{1} r_1+\alpha_{2} r_2+\ldots+\alpha_{n} r_n\right)\frac{\partial}{\partial r_i}.$$  

Notice that at least one $v_i$ is nonzero since otherwise we obtain $D_1=0$. Consequently, $\a_i\neq 0$ for some $i$. Suppose that $\a_1\neq 0$. 
 Denote by $J$ the ideal of $R_n$ generated by $r_2,\ldots, r_n$. Note that $D_1^{(r)}(J)\subseteq J$. Consider the factor algebra $R_n/J\cong K\left\langle r_1\right\rangle$. Then $D_1^{(r)}$ induces the derivation $D_0=\alpha_{1} r_1 r_1\frac{\partial}{\partial r_1}$ on $K\left\langle r_1\right\rangle$.  Obviously, $D_0$ is not locally nilpotent. Consequently, $D_1^{(r)}$ is not locally nilpotent. 
 $\Box$

\begin{lemma}\label{NLNlem10} If $g=f$ and $v_1=0$, then $D_1^{(r)}$ is not locally nilpotent. 
\end{lemma}
\Proof  If $g=f$ and $v_1=0$,  then the equality  \eqref{NLN7} implies
\begin{gather} \label{NLN15}
r_1 u_j=u_jr_1+r_1 v_j \ (\text{if} \ i=1, j\neq 1),
\end{gather}
\begin{gather} \label{NLN16}
r_1 u_1=u_1 r_1 \ (\text{if} \ i=j=1).
\end{gather}

Recall that, as mentioned in the proof of Theorem \ref{tu1}, $D(l_i)=\mathrm{ad}_{l_i} (w)$ for some  $w\in I_n$ and 
$\mathrm{Lm}(D(l_i))=f=\mathrm{Lm}(w)\neq 0$ for all $1 \leq i \leq n$. Then the condition $g=f$ implies also that $u_i=\mathrm{ad}_{l_i}(\mathrm{Lc}(w))$. Set $w_0=\mathrm{Lc}(w)$. Then $u_i=\mathrm{ad}_{l_i}(w_0)$. 

Since $w\in I_n$ it follows that $\deg(w_0)\geq 1$ and $\deg (u_i) =t\geq 2$.
Then $v_i\in (U_n)_t$ since $D$ is homogenous.

By Bergman’s Theorem \cite{Bergman} on centralizers, the equation \eqref{NLN16} implies $u_1\in K\left\langle  r_1\right\rangle$. Therefore,  
$$u_1=\beta r_1^t=-\frac{\beta}{t-1} \mathrm{ad}_{l_1}(r_1^{t-1}), \  \beta\in K^*.$$ 
This implies $w_0=-\frac{\beta}{t-1} r_1^{t-1}$ since the kernel of  $\mathrm{ad}_{l_1}$ in $R_n$ is zero. 
Consequently, 
$$u_i=-\frac{\beta}{t-1} \mathrm{ad}_{l_i}( r_1^{t-1}) \text{ for all } 1\leq i \leq n.$$
By substituting $u_j$, $j\neq 1$, into the equality \eqref{NLN15}, we obtain 
$$-\frac{\beta}{t-1} r_1 \mathrm{ad}_{l_j}( r_1^{t-1})=-\frac{\beta}{t-1} \mathrm{ad}_{l_j}( r_1^{t-1})r_1+r_1 v_j.$$
This can be rewritten as  
$$ r_1 (r_1r_jr_1^{t-2}+r_1^2r_jr_1^{t-3}+\ldots+r_1^{t-2}r_jr_1+r_1^{t-1}r_j)$$
$$= (r_1r_jr_1^{t-2}+r_1^2r_jr_1^{t-3}+\ldots+r_1^{t-2}r_jr_1+r_1^{t-1}r_j)r_1+\frac{t-1}{\beta}r_1 v_j.$$
From this we obtain
$$ r_1^{t} r_j=  r_1r_j r_k^{t-1}+\frac{t-1}{\beta}r_1 v_j.$$
Hence
$$v_j=\frac{\beta}{t-1}(r_1^{t-1} r_j-r_j r_1^{t-1})\text{ for all } 2\leq j \leq n.$$
Therefore, 
 $$D_1^{(r)}=\sum^{n}_{i=2}\frac{\beta}{t-1}[r_1^{t-1}, r_i]\frac{\partial}{\partial r_i}, \ t>1.$$  
Notice that 
$$ (D_1^{(r)})^{d}(r_2)=\frac{\beta^d}{(t-1)^d} \mathrm{ad}_{r_1^{t-1}}(r_2)\neq 0$$ 
for all $d\geq 0$, that is, $D_1^{(r)}$ is not locally nilpotent. $\Box$

\begin{lemma}\label{NLNlem9} If  $g=f$ and $v_i\neq 0$ for all $i$, then $D_1^{(r)}$ is not  locally nilpotent.
\end{lemma}
\Proof If  $g=f$ and $v_i\neq 0$ for all $i$, then the equality \eqref{NLN7} implies 
\begin{gather} \label{NLN8}
-\mathrm{ad}_{l_j}(v_i)-v_ir_j+r_i u_{j}=u_j r_i+r_i v_{j}.
\end{gather}
 In the proof of Lemma \ref{NLNlem10} under the condition $g=f$ was shown that there exists $0\neq w_0\in R_n\cup I_n$ such that $u_i=\mathrm{ad}_{l_i}(w_0)$ for all $i$. Then $\deg (u_i) =t\geq 2$. Let 
$$w_0=\sum^{n}_{s=1}r_sb_{s}=\sum^{n}_{t=1}c_{t}r_t, \ \ 
v_{i}=\sum^{n}_{s=1} \sum^{n}_{t=1}r_s v^i_{st} r_t,$$
where $b_{s}, c_{t}, v^i_{st} \in R_n$ for all $s,t$.

 By substituting $v_i$ into the equality \eqref{NLN8}, we obtain 
$$-\sum^{n}_{s=1} \sum^{n}_{t=1}\mathrm{ad}_{l_j}( r_sv^i_{st}) r_t+\sum^{n}_{s=1} \sum^{n}_{t=1}r_s v^i_{st} r_tr_j-\sum^{n}_{s=1} \sum^{n}_{t=1}r_s v^i_{st} r_tr_j +r_i u_{j}=u_j r_i+r_i v_{j}.$$
Therefore 
$$-\sum^{n}_{s=1} \sum^{n}_{t=1}\mathrm{ad}_{l_j}(r_s v^i_{st}) r_t +r_i u_{j}=u_j r_i+r_i v_{j}.$$
Hence $\mathrm{ad}_{l_j}(r_s v^i_{st})=0$, and, consequently, $v^i_{st}=0$ if $s,t\neq i$. Thus
 $$v_{i}= \sum^{n}_{t=1}r_i v^i_{it} r_t+\sum^{n}_{s=1, s\neq i}r_s v^i_{si} r_i.$$

By substituting $v_i$, $u_j$ into the equality \eqref{NLN8}, we obtain 
\begin{equation}\label{eq_12}
-\sum^{n}_{t=1} \mathrm{ad}_{l_j}(r_iv^i_{it}) r_t-\sum^{n}_{s=1, s\neq i} \mathrm{ad}_{l_j}(r_sv^i_{si}) r_i +r_i u_{j}=\mathrm{ad}_{l_j}(\sum^{n}_{s=1}r_sb_{s}) r_i+r_i v_{j}. 
\end{equation}
Hence
$$ \sum^{n}_{s=1, s\neq i}r_sb_{s} =-\sum^{n}_{s=1, s\neq i}r_sv^i_{si},$$
and, consequently, $b_{s}=v^i_{si}$ if $s\neq i$. 
Therefore,  we can write $v_i$ as
$$ v_i=\sum^{n}_{t=1}r_i v^i_{it} r_t+r_ib_ir_i-w_0r_i$$
since $\sum^{n}_{s=1, s\neq i}r_sb_{s}=w_0-r_ib_i$.

The equality \eqref{eq_12} implies 
$$ -\sum^{n}_{t=1}\mathrm{ad}_{l_j}(r_i  v^i_{it}) r_t +r_i u_{j}=\mathrm{ad}_{l_j}(r_ib_i)r_i+r_i v_{j}.$$
The case $j=i$  can be written as
$$\sum^{n}_{t=1}r_i r_i v^i_{it} r_t-\sum^{n}_{t=1}r_i \mathrm{ad}_{l_i}( v^i_{it}) r_t +r_i ad_{l_i}(\sum^{n}_{t=1}c_{t}r_t)$$
$$=-r_ir_ib_ir_i+r_i\mathrm{ad}_{l_i}(b_i)r_i+r_i (\sum^{n}_{t=1}r_i v^i_{it} r_t+r_ib_ir_i-w_0r_i). $$
From this we obtain 
\begin{equation}\label{eq_13}
-r_i\sum^{n}_{t=1}\mathrm{ad}_{l_i}(  v^i_{it}) r_t +r_i \sum^{n}_{t=1}\mathrm{ad}_{l_i}\left(c_{t}\right)r_t-r_i \sum^{n}_{t=1}c_{t}r_tr_i=r_i\mathrm{ad}_{l_i}(b_i)r_i-r_i w_0 r_i.
\end{equation}
Hence $$ \sum^{n}_{t=1, t\neq i}c_{t}r_t=\sum^{n}_{t=1, t\neq i} v^i_{it} r_t,$$
and, consequently, $c_{t}=v^i_{it}$ if $t\neq i$. 
Therefore,  we can write $v_i$ in the form 
$$ v_i=r_i v^i_{ii} r_i-r_ic_ir_i+r_iw_0+r_ib_ir_i-w_0r_i=r_i(v^i_{ii}-c_i+b_i)+r_iw_0-w_0r_i$$
since $ \sum^{n}_{t=1, t\neq i}c_{t}r_t=w_0-c_ir_i$. Note that $v^i_{ii}-c_i+b_i \in (U_n)_k$ where $k=\deg(w_0)\geq 1$. 

The equality \eqref{eq_13} implies 
$$ -r_i\mathrm{ad}_{l_i}(  v^i_{ii}) r_i +r_i \mathrm{ad}_{l_i}(c_i)r_i-r_iw_0r_i=r_i\mathrm{ad}_{l_i}(b_i)r_i-r_i w_0 r_i,$$
which gives 
$$r_i\mathrm{ad}_{l_i}(-v^i_{ii}+c_i-b_i)r_i=0. $$
It follows that $-v^i_{ii}+c_i-b_i= 0$. Consequently, 
$$v_{i}= r_iw_0-w_0r_i=-\mathrm{ad}_{w_0}(r_i).$$ 
It means $ D_1^{(r)}=-\mathrm{ad}_{w_0}$, and is not locally nilpotent. $\Box$

\begin{pr}\label{th1}
A locally nilpotent derivation $D_0=g(l_n)\frac{\partial}{\partial l_1}$ of $L_n$, $n\geq 2$,  can be uniquely extended to a locally nilpotent derivation $D=g(l_n)\frac{\partial}{\partial l_1}+g'(l_n)r_n\frac{\partial}{\partial r_1} $ of $U_n$. 
\end{pr}

\Proof First we show that $D$ is indeed a derivation. It is sufficient to check that $D$ preserves the relations \eqref{s1} and \eqref{s2}. Applying $D$ to the relation \eqref{s2}, we obtain
\bee\label{equ4}
D(r_i)l_j+r_iD(l_j)=D(l_j)r_i+l_jD(r_i)+D(r_i)r_j+r_iD(r_j).
\eee
This is true for all $i,j\geq 2$ since $D(l_i)=D(r_i)=0$ for all $i$. For the remaining values of $i,j$  \eqref{equ4}  implies 
$$ g'(l_n)r_nl_1+r_1g(l_n)=g(l_n)r_1+l_1g'(l_n)r_n+g'(l_n)r_nr_1+r_1g'(l_n)r_n,  $$
$$ g'(l_n)r_nl_j=l_jg'(l_n)r_n+g'(l_n)r_nr_1, $$
$$r_ig(l_n)=g(l_n)r_i+r_ig'(l_n)r_n.$$
By relation \eqref{s2} and Corollary  \ref{c_1}, these equalities hold.
Since  $D^2(l_1)=D^2(r_1)=0$ it follows that $D$ is locally nilpotent.

Let $D$ be arbitrary locally nilpotent derivation of $U_n$ that induces  $D_0$ on $U_n/I_n\cong L_n$. Recall that every derivation of $U_n$ induces a derivation of $L_n$ by Lemma \ref{TRlem1}. Then 

$$D=\left(g(l_n)+a_1\right)\frac{\partial}{\partial l_1}+\sum^{n}_{k=2}a_k\frac{\partial}{\partial l_k}+\left(g'(l_n)r_n+b_1\right)\frac{\partial}{\partial r_1}+\sum^{n}_{k=2}b_k\frac{\partial}{\partial r_k},$$
where $a_i,b_i \in I_n$ for all $i$.  
We show that $a_i=b_i=0$ for all $i$. 

Assume the contrary. Consider the $w=(1,\ldots, 1, 0)$-degree function and the highest $w$-homogeneous part $\widehat{D}$ of $D$. If $\mathrm{wdeg}(D)\geq 0$, then 
$$\widehat{D}= \sum^{n}_{k=1}c_k\frac{\partial}{\partial l_k}+\sum^{n}_{k=1}d_k\frac{\partial}{\partial r_k},$$
 where $c_i,d_i \in I_n$ for all $i$. By Lemma \ref{GRpr2}, the derivation $\widehat{D}$ is locally nilpotent. By Theorem \ref{tu1}, $\widehat{D}=0$ since $D_1(U_n)\subseteq I_n$. 

Consequently, $\mathrm{wdeg}(D)=-1$. It is possible only if $a_n=b_n=0$ and  $a_i,b_i$ belong to the subalgebra $\langle l_n,r_n\rangle$ generated by 
$l_n, r_n$ for all $1\leq i\leq n$.

Applying $D$ to the relation \eqref{s1} we obtain $\mathrm{ad}_{l_i}(D(l_j))=\mathrm{ad}_{l_j}(D(l_i)).$  It follows that $\mathrm{ad}_{l_i}(a_j)=\mathrm{ad}_{l_j}(a_i).$  
By Lemma \ref{NLNlem2}, there exists $h\in I_n$ such that $a_i=\mathrm{ad}_{l_i}(h)$ for all $i$. If $h\neq 0$, then $a_1=\mathrm{ad}_{l_1}(h)\neq 0$ cannot belong to the subalgebra $\langle l_n,r_n\rangle$. Consequently, $w=0$ and $a_1=\ldots=a_n=0$. Then 
  $$D(l_1)=g(l_n) \text{  and  } D(l_k)=0, \ 2\leq k\leq n.$$

Applying $D$ to the defining relation $r_1l_k=l_kr_1+r_1r_k$, $2\leq k\leq n$, we obtain
$$(g'(l_n)r_n+b_1)l_k=l_k(g'(l_n)r_n+b_1)+(g'(l_n)r_n+b_1)r_k+r_1b_k. $$
This can be written as  
$$g'(l_n)l_kr_n+g'(l_n)r_nr_k+l_kb_1-ad_{l_k}(b_1)=g'(l_n)l_kr_n+l_kb_1+g'(l_n)r_nr_k+b_1r_k+r_1b_k. $$
Consequently, 
\begin{equation}{\label{eq_17}}
-\mathrm{ad}_{l_k}(b_1)=b_1r_k+r_1b_k.
\end{equation}
It follows that $b_k=0$ for all $2\leq k\leq n$ since $b_1\in\langle l_n,r_n\rangle$. 
We can represent $b_1$ as $b_1=qr_n$ since $b_1\in I_n$. Then the  equality \eqref{eq_17}  implies $$-\mathrm{a}d_{l_k}(q)r_n+qr_nr_k=qr_nr_k.$$
Hence $\mathrm{ad}_{l_k}(q)=0$. It follows that $q\in K[l_n]$ since $b_1\in\langle l_n,r_n\rangle$. Consequently,  
$$ D(r_1)=g'(l_n)r_n+q(l_n)r_n \text{  and  } D(r_k)=0, \ 2\leq k\leq n.$$ 

Applying $D$ to  the defining relation $ r_kl_1=l_1r_k+r_k r_1$, $2\leq k\leq n$, we obtain 
$$ r_kg(l_n)=g(l_n)r_k+ r_k(g'(l_n)r_n+q(l_n) r_n).$$ 
By Corollary \ref{c_1}, we can write this equality as
$$g(l_n)r_k+r_kg'(l_n)r_n=g(l_n)r_k+ r_kg'(l_n)r_n+r_kq(l_n) r_n. $$
From this we obtain $ r_k q(l_n)r_n=0$. Hence $q(l_n)=0$. Consequently,  $D(r_1)=g'(l_n)r_n$. $\Box$

\section{Automorphisms of $U_n$, $n\geq2$}
\hspace*{\parindent}

The main result of our paper is the following theorem. 
\begin{theorem}\label{th2}
 $\mathrm{Aut}(U_n)\cong \mathrm{Aut}(L_n), \ n\geq 2 .$
\end{theorem}
\Proof Let  $\theta$ be an arbitrary automorphism of $U_n$. By Lemma \ref{ELlem1}, $\theta(I_n)\subseteq I_n$. Thus, $\theta$ induces an automorphism $\varphi$ of $U_n/I_n\cong L_n$.  By  Proposition \ref{up1}, $\varphi$ can be lifted to an automorphism $\varphi^*$ of $U_n$. 
 Let  $\theta'=\theta(\varphi^*)^{-1}$. We are going to show that $\theta'=id$, and, consequently, $\theta=\varphi^*$. This proves that the embedding $\Phi$ in the Proposition 
\ref{up1} is an isomorphism of groups.

Obviously,  
$$\theta'(l_i)=l_i+a_i,\  \theta'(r_i)=r_i+b_i,$$
where $a_i,b_i\in I_n$. Note that, by Lemma \ref{affaut}, $a_i,b_i\in \oplus_{ m\geq 2} (U_n)_m$. 

Since  $l_i+a_i, r_i+b_i$ are the generators of $U_n$ we can define a derivation $D_h$ given by 
$$D_h(l_1+a_1)=h(l_n+a_n), \ D_h(l_k+a_k)=0,$$
$$D_h(r_1+b_1)=h'(l_n+a_n)(r_n+b_n), \ D_h(r_k+b_k)=0,$$
where $h\in k[x]$, $2\leq k \leq n$. 
This derivation is locally nilpotent. We have
$$D_h(l_1)=h(l_n+a_n)-D_h(a_1)=h(l_n)+a'_1,$$
$$D_h(l_k)=-D_h(a_k)=a'_k,$$
$$D_h(r_1)=h(l_n+a_n)(r_n+b_n)-D_h(b_1)=h'(l_n)r_n+b'_1,$$
$$D_h(r_k)=-D_h(b_k)=b'_k,$$
where $a'_i,b'_i\in I_n$, since the ideal $I_n$ is invariant under every derivation. Proposition \ref{th1} gives that $a'_i=b'_i=0$ since $D_h$ is locally nilpotent. Therefore 
$$D_h=h(l_n)\frac{\partial}{\partial l_1}+h'(l_n)r_n\frac{\partial}{\partial r_1}$$
and
$$D_h(a_1)=h(l_n+a_n)-h(l_n), \ D_h(a_k)=0,$$
$$D_h(b_1)=h'(l_n+a_n)(r_n+b_n)-h'(l_n)r_n, \ D_h(b_k)=0.$$
If $\deg(h)=m$, then $\deg (D_h)=m-1$ and it follows from Lemma \ref{GRlem1} that $\deg (D_h(g))\leq \deg (g) +m-1$ for every $g\in U_n$.

Consider $h=x$. Then
$$D_h(a_1)=a_n, \ D_h(a_k)=0, \ D_h(b_1)=b_n, \ D_h(b_k)=0.$$
Therefore $\deg (D_h(a_1))\leq\deg (a_1)$ and $\deg (D_h(b_1))\leq\deg (b_1).$ Hence $\deg (a_n)\leq\deg (a_1)$ and $\deg (b_n)\leq\deg (b_1).$ We also obtain  
$\deg (a_1) \leq \deg (a_n)$ and $\deg (b_1) \leq \deg (b_n) $ if,  in the definition of $D_h$  we exchange $l_1$ with $l_n$ and $r_1$ with $r_n$. Consequently, 
$$\deg (a_1)=\deg (a_n), \ \ \deg (b_1)=\deg (b_n).$$
We also obtain  
\begin{equation}{\label{E1}}
\deg (a_k)=\deg (a_1), \ \ \deg (b_k)=\deg (b_1), \ 2 \leq k\leq n,
\end{equation}
 if, in the definition of $D_h$  we exchange $l_n$ with $l_k$ and $r_n$ with $r_k$.

We show that $a_i=0$ for all $i$. Assume that this is not the case. Set $h=x^2$. Then 
$$D_h(a_1)=(l_n+a_n)(l_n+a_n)-l^2_n=a_nl_n+l_na_n+a_n a_n$$
and
$$\deg (a_nl_n+l_na_n+a_n a_n)\leq \deg a_1 +1.$$
Therefore $\deg (a_1)+1\geq 2\deg (a_n)$ since $\deg (a_n)\geq 2$. By \eqref{E1}, $\deg (a_1)+1\geq 2\deg (a_1)$. Hence $\deg (a_1)\leq 1$. This contradiction gives $a_i=0$ for all $i$.

Now we show that $b_i=0$ for all $i$. Assume that this is not the case.  Applying the automorphism $\theta'$ to the relation $r_1l_1=l_1r_1+r_1r_1$, we obtain
$$(r_1+b_1)l_1=l_1(r_1+b_1)+(r_1+b_1)(r_1+b_1).$$
By \eqref{s2}, we have
$$l_1 r_1+r_1r_1 +b_1 l_1=l_1 r_1+l_1 b_1+r_1 r_1+r_1 b_1+b_1 r_1+b_1b_1.$$
From this we obtain
$$-\mathrm{ad}_{l_1}(b_1)=r_1 b_1+b_1 r_1+b_1b_1.$$
This is impossible since since $\deg (b_1)\geq 2$. This contradiction gives $b_i=0$ for all $i$. Consequently, $\theta'=id$. $\Box$

\section{Derivations and automorphisms of $U_1$}
\hspace*{\parindent}

This section is devoted to describing the locally nilpotent derivations and automorphisms of $U_1$. Recall that $L_1=K[l_1]$ and $R_1=K[r_1]$.

\begin{lemma}\label{U12} Let $D=(u+\gamma)\frac{\partial}{\partial l_1}+v \frac{\partial}{\partial r_1}$, where $u,v\in I_1$, $\gamma\in K$, be a non-zero derivation of $U_1$.  Let  $\mathrm{Lm}(u)=l_1^m$ and $\mathrm{Lm}(v)=l^s_1$.  Then 
\begin{enumerate}
\item  $u\in R_1$ and $v=0$ if  $m-1>s$;
\item  $\mathrm{Lc}(v)=\alpha r_1^2$, $\alpha\in K^*$, if $m-1<s$;
\item  $\widehat{\mathrm{Lc}(u)}=\frac{-p+2}{s+1}\beta r_1^{p-1}$, $\widehat{\mathrm{Lc}(v)}=\beta r_1^p$,   $p>2$, $\beta \in K^*$, if $m-1=s$.
\end{enumerate}
\end{lemma}
\Proof The derivation $\frac{\partial}{\partial l_1}$ of $L_1$ can be continued to a derivation of $U_1$ by $\frac{\partial}{\partial l_1}(r_1)=0$. Denote this derivation of $U_1$ by $\d_1$. Formally, $\d_1=\frac{\partial}{\partial l_1}$ in the form (\ref{u3}). 

By Corollary \ref{c_1}, we have
\bes r_1l_1^k r_1^t=l_1^k r_1^{t+1}+kr_1l_1^{k-1} r_1^{t+1}=l_1^k r_1^{t+1}+r_1\partial_1(l_1^k r_1^t)r_1.
\ees
Then
\bee\label{equ5}
r_1w=wr_1+r_1 \partial_1(w)r_1 
\eee
for all   $w\in U_1$.

Applying $D$ to the relation \eqref{s2}, we obtain $vl_1+r_1u=ur_1+l_1v+vr_1+r_1v.$  By \eqref{equ5}, we can write this as
$$l_1v-\mathrm{ad}_{l_1}(v)+ur_1+r_1\partial_1(u)r_1=ur_1+l_1v+vr_1+vr_1+r_1\partial_1(v)r_1.$$
From this we obtain
\bes
-\mathrm{ad}_{l_1}(v)-2vr_1-r_1\partial_1(v)r_1+r_1\partial_1(u)r_1=0.
\ees

If $\mathrm{Lc}(u)=u_m$ and $\mathrm{Lc}(v)=u_s$, then, by \eqref{equ5}, the last equality implies 
\begin{equation} \label{U15}
l_1^s(-\mathrm{ad}_{l_1}(v_s)-2v_sr_1)+f +ml_1^{m-1}r_1r_1u_m +g=0,
\end{equation}
where $\mathrm{Lm}(f)<l_1^s$, $\mathrm{Lm}(g)<l_1^{m-1}$.

(1) If $m-1>s$, then,  by Lemma \ref{LTlem3}, \eqref{U15} implies $mr_1r_1u_m=0$. Hence $u_m=0$ if $m>0$. Consequently, $m=0$, $ u\in  R_1$ and  $v=0$ since $m-1>s$.

(2) If $m-1<s$, then \eqref{U15} implies $- ad_{l_1}(v_s)=2 r_1 v_s$. By Lemma \ref{U_14}, $\mathrm{Lc}(v)=v_s=\alpha r_1^2$, $\alpha\in K^*$.

(3)  If $m-1=s$, then, by Lemma \ref{LTlem3}, \eqref{U15} implies 
$$- ad_{l_1}(v_s)-2 r_1 v_s+m r_1r_1u_m=0.$$
Let $$u_{m}=\alpha r_1^t+a, \  \ v_s=\beta r_1^p +b,$$
 where $\alpha, \beta \in K^*$, $a,b \in R_1\cap I_1$, $\deg(a)<t$,  $\deg(b)<p$.
Note that $t,p\geq 1$ since $u_m, v_s \in R_1\cap I_1$.  Substituting them into the last equation, we obtain 
$$ (p-2)\beta  r_1^{p+1}- ad_{l_1}(b)-2r_1b+m\alpha r_1^{t+2}+mr_1^2a=0.$$
It follows that $p\neq 2$ and $t+2=p+1$, and, consequently, $p=t+1>2$. Hence  
$ \alpha=\frac{-p+2}{m}\beta=\frac{-p+2}{s+1}\beta.$ $\Box$

\begin{theorem}\label{U10} Let $D$ be a locally nilpotent derivation of $U_1$. Then $D(l_1)\in R_1 $ and $D(r_1)=0$.
\end{theorem}
\Proof Let $D$ be a non-zero locally nilpotent derivation of $U_1$. By Lemma \ref{TRlem1}, $D(I_1)\subseteq I_1$. So $D$ induces a locally nilpotent derivation $D_0$ on $U_1/I_1\cong  L_1$. Hence $D_0=\gamma \frac{\partial }{\partial l_1}$, $\gamma\in K$, and $D=(u+\gamma)\frac{\partial}{\partial l_1}+v \frac{\partial}{\partial r_1}$, where $u,v\in I_1$. Let  $\mathrm{Lm}(u)=l_1^m$ and $\mathrm{Lm}(v)=l^s_1$.

a)  If  $m-1>s$, then,  by Lemma \ref{U12},   $u\in R_1$ and $v=0$. Obviously, in this case $D$ is locally nilpotent.

b)  If $m-1<s$,  then,  by Lemma \ref{U12},  $\mathrm{Lc}(v)=\alpha r_1^2$, $\alpha\in K^*$.  It is clear that $\mathrm{Lc}(D)^{(r)}=\alpha r_1^2$ and $\mathrm{Lc}(D)^{(r)}$ is not locally nilpotent derivation of $R_1$. By Lemma \ref{NLNpr2},  $D$ is not locally nilpotent. 

c) If  $m-1=s$, then,  by Lemma \ref{U12},    $$D(l_1)=\frac{-p+2}{s+1}\beta l_1^{s+1}r_1^{p-1}+l_1^{s+1}a+c, \ \  D(r_1)=\beta l_1^{s}r_1^p+l_1^{s}b+d, \ \ p>2, \ \beta\in K^*, $$
where $a, b \in  R_1 \cap I_1$, $\deg(a)<p-1$, $\deg(b)<p$, $c, d\in U_1$,   $\mathrm{Lm}(c)<l_1^{s+1}$, $\mathrm{Lm}(d)<l_1^{s}$.

By induction on $k\geq 1$, we want to prove that 
\begin{equation} \label{eq_18}
 D^k(r_1)= \gamma_k l_1^{ks}r_1^{k(p-1)+1}+l_1^{ks}h_k+g_k,
\end{equation}
where $\gamma_k= \frac{\beta^k}{(s+1)^{k-1}} \prod^{k}_{i=1} (ks+(k-1)p-(k-2))$, $h_k\in R_1 \cap I_1$, $\deg(h_k)<k(p-1)+1$, $g_k \in U_1$,  $\mathrm{Lm}(g_k)<l_1^{ks}$. 

By the induction proposition, Lemma \ref{LTlem3}, and Lemma \ref{GRlem1}, we have
$$D^{k+1}(r_1)=\gamma_k ks l_1^{ks-1}\left(\frac{-p+2}{s+1}\beta l^{s+1}r_1^{p-1}\right)r_1^{k(p-1)+1}$$
$$+\gamma_k (k(p-1)+1)l_1^{ks}(\beta l^{s}r_1^p) r^{(k+1)(p-1)}+l_1^{(k+1)s}h_{k+1}+g_{k+1} $$
$$=\gamma_{k}\frac{\beta}{s+1}((k+1)s+kp-(k-1))l_1^{(k+1)s}r_1^{(k+1)(p-1)+1}+l_1^{(k+1)s}h_{k+1}+g_{k+1}$$
$$=\gamma_{k+1}l_1^{(k+1)s}r_1^{(k+1)(p-1)+1}+l_1^{(k+1)s}h_{k+1}+g_{k+1},$$
where $h_{k+1}\in R_1 \cap I_1$, $\deg(h_{k+1})<(k+1)(p-1)+1$, $g_{k+1} \in U_1$,   $\mathrm{Lm}(g_{k+1})<l_1^{(k+1)s}$. Thus,  the equation \eqref{eq_18} holds for any $k\geq 1$.
It is clear that $\gamma_k\neq 0$ for all $k\geq 1$. Consequently, $D$ is not locally nilpotent. $\Box$

\begin{theorem}\label{U11} Let  $\phi$ be an automorphism of $U_1$. Then
 $$ \phi(l_1)=\alpha l_1+h(r_1),  \ \varphi(r_1)=\alpha r_1, \text{ where } \alpha\in K^\ast, \ \ h(r_1)\in R_1.$$
\end{theorem}
\Proof   Let  $\phi$ be an arbitrary automorphism of $U_1$.  By Lemma \ref{ELlem1}, $\phi(I_1)\subseteq I_1$. So $\phi$ induces an automorphism $\psi$ of $U_1/I_1\cong L_1$. This means that $\psi(l_1)=\alpha l_1+\beta$, $ \alpha, \beta\in K$, $\alpha\neq 0$. By Proposition  \ref{up1}, $\psi$ can be lifted to an automorphism $\psi^*$ of $U_1$.  Set $\phi'=\phi(\psi^*)^{-1}$. Then $\phi'$ induces the identity automorphism of $L_1$. Therefore, 
$$\phi'(l_1)=l_1+a, \ \phi'(r_1)=r_1+b, \text{ where } a,b\in I_1.$$

Since  $l_1+a, r_1+b$ are  generators of $U_1$ we can define a locally nilpotent derivation $D_h$ as
$$D_h(l_1+a)=h(r_1+b), \ D_h(r_1+b)=0, \ \text{where} \ h\in K[x].$$ 
We have
$$D_h(l_1)=h(r_1+b)-D_h(a)=c, \ \ D_h(r_1)=-D_h(b)=d,$$
where $c,d\in I_1$, since the ideal $I_1$ is invariant under derivations. By Theorem \ref{U10}, if $D_h$ is locally nilpotent, then $c=c(r_1)\in R_1$ and $d=0$. Therefore, 
$$D_h=c(r_1)\frac{\partial}{\partial l_1}, \ \ D_h(a)=h(r_1+b)-c(r_1), \ \ D_h(b)=0.$$
Hence $b \in R_1$. This means that $\phi(r_1)=r_1+b$ determines an automorphism of the free associative algebra $R_1$. Consequently,  $b=0$.

Then $D_h(a)=h(r_1)-c(r_1)\in R_1$. Hence $a=l_1a_1+a_2$, $a_1, a_2 \in R_1 \cap I_1$.
Applying $\phi$ to the relation \eqref{s2} , we obtain 
$$ r_1l_1+r_1l_1a_1+r_1a_2=l_1r_1+l_1a_1r_1+a_2r_1+ r_1r_1.$$
From this we obtain $r_1r_1a_1=0$ and $a_1=0$. Consequently, $a \in R_1.$ $\Box$

\section*{Acknowledgments}

This work was supported by the grant AP23487886 of the Ministry of Science and Higher Education of the Republic of Kazakhstan.

\end{document}